\documentclass[12pt,twoside]{amsart}
\usepackage{amssymb}
\usepackage{amscd}

\title{On Kawamata's theorem}
\author{Osamu Fujino} 
\date{2009/6/3}
\address{Department of Mathematics, Faculty of 
Science, Kyoto University, Kyoto 606-8502 Japan}
\email{fujino@math.kyoto-u.ac.jp}

\newcommand{\rank}[0]{{\operatorname{rank}}}

\newcommand{\h}[0]{{\operatorname{hom}}}
\newcommand{\KC}[0]{{\operatorname{KC}}}
\newtheorem{thm}{Theorem}[section]
\newtheorem{lem}[thm]{Lemma}
\newtheorem{cor}[thm]{Corollary}

\theoremstyle{definition}
\newtheorem{defn}[thm]{Definition}
\newtheorem{rem}[thm]{Remark}
\newtheorem*{ack}{Acknowledgments}       
 
\newtheorem{say}[thm]{}
\begin{document}
\bibliographystyle{amsalpha+}

\begin{abstract}
We give an alternate proof of 
the main theorem of 
Kawamata's 
paper:~Pluricanonical systems on minimal algebraic varieties. 
Our proof also works for varieties in class $\mathcal C$. 
We note that our proof is completely different from Kawamata's. 
\end{abstract}

\maketitle 
\tableofcontents

\section{Introduction} 
One of the main purposes of this paper is to cut a chain of 
troubles caused by \cite[Theorem 4.3]{kawamata}. 
We give an alternate proof of the following famous theorem, 
which we call {\em{Kawamata's theorem}} 
in this paper. This theorem is indispensable for 
the abundance conjecture. 
\begin{thm}[{cf.~\cite[Theorem 6-1-11]{kmm}}]\label{thm2}
Let $(X,B)$ be a klt pair and $\pi:X \to S$ a proper 
surjective morphism of normal varieties. Assume the following 
conditions: 
\begin{itemize}
\item[(a)] $H$ is a $\pi$-nef $\mathbb Q$-Cartier divisor on $X$, 
\item[(b)] $H-(K_X+B)$ is $\pi$-nef and $\pi$-abundant, and 
\item[(c)] $\kappa (X_\eta, (aH-(K_X+B))_\eta)\geq 0$ 
and $\nu(X_\eta, (aH-(K_X+B))_\eta)=\nu(X_\eta, (H-(K_X+B))_\eta)$ for 
some $a\in \mathbb Q$ with $a>1$, 
where $\eta$ is the generic point of $S$. 
\end{itemize} 
Then $H$ is $\pi$-semi-ample. 
\end{thm}
It was first proved in \cite{kawamata} on the 
assumption that $S$ is a point. 
Kawamata's proof heavily depends on a very technical 
generalization of Koll\'ar's injectivity theorem on {\em{generalized 
normal crossing varieties}} (see \cite[Section 4]{kawamata}). 
Once we adopt this difficult injectivity theorem, X-method 
works and the proof is essentially the same as 
the one of the Kawamata--Shokurov base point free theorem. 
Unfortunately, there is an ambiguity in the proof of 
\cite[Theorem 4.3]{kawamata} 
(see \cite[Remark 3.10.3]{fujino} and \ref{511} below). 
Thus, our proof is the first rigorous proof of Kawamata's 
theorem. 
It is completely different from Kawamata's. 
His proof relies on the theory of mixed Hodge structures for 
reducible varieties. Our proof grew out from the theory of 
variation of Hodge structures, especially, Deligne's 
canonical extensions of Hodge bundles. 
We note that our method saves Kawamata's 
theorem but does not recover the results in \cite[Section 4]
{kawamata}. 
They are completely generalized 
in \cite[Chapter 2]{fuji-vani} 
for {\em{embedded simple normal crossing pairs}}. 
However, 
\cite{fuji-vani} does not recover \cite[Theorem 4.3]{kawamata}. 
Compare the arguments in \cite[Chapter 2]{fuji-vani} with 
Kawamata's ones. 
The reader can find a slight generalization of Kawamata's 
theorem and some other applications of our methods in 
\cite{base} and \cite{fujino5}. 

We summarize the contents of this paper. 
In Section \ref{sec2}, we will give an alternate 
proof of 
Kawamata's theorem. By using Ambro's formula, 
we will reduce Kawamata's theorem to 
a reformulated version of the Kawamata--Shokurov base 
point free theorem. 
Section \ref{formula} is an appendix, where 
we will quickly review Ambro's formula for the 
reader's convenience. In Section \ref{sec4}, 
we will prove Kawamata's theorem for 
varieties in class $\mathcal C$, which is 
\cite[Theorem 5.5]{nakayama1}. 
We separate this section from Section \ref{sec2} in order not to 
make needless confusion. 
In the final section:~Section \ref{sec5}, we will make 
some comments on topics related to Kawamata's 
theorem for the coming generation. 

\begin{ack}
I was partially supported by The Sumitomo Foundation, The Inamori 
Foundation, and 
by the Grant-in-Aid for Young Scientists (A) 
$\sharp$20684001 
from JSPS. 
I would like to thank Professors Yujiro Kawamata and 
Noboru Nakayama for 
answering my questions. 
I thank Professors Daisuke Matsushita and 
Shigefumi Mori for 
encouraging me during the preparation of this 
paper. 
\end{ack}

We will work over an algebraically closed field $k$ of 
characteristic zero throughout 
this paper. We adopt the language of {\em{b-divisors}} and 
use the standard notation of the log minimal model program. 
See, for example, \cite{corti}. 

\section{Proof of Kawamata's theorem}\label{sec2}

The following theorem is a reformulation of the Kawamata--Shokurov 
base point free theorem. 
The original proof works without any changes (cf.~\cite[Theorem 3-1-1]{kmm}). 

\begin{thm}[Base point free theorem]\label{thm1}
Let $(X,B)$ be a sub klt pair, let $\pi:X \to S$ be a proper 
surjective morphism of normal varieties, and $D$ a $\pi$-nef 
Cartier divisor on $X$. Assume the following conditions: 
\begin{itemize}
\item[(1)] $rD-(K_X+B)$ is nef and big over $S$ for some positive integer 
$r$, and 
\item[(2)] $\pi_*\mathcal O_X(\ulcorner \mathbf 
A(X,B)\urcorner +j\overline D)\subseteq 
\pi_*\mathcal O_X(jD)$ for every positive integer $j$, 
where $\mathbf A(X, B)$ is the discrepancy 
$\mathbb Q$-b-divisor and 
$\overline D$ is the Cartier 
closure of $D$ $($see \cite[Example 2.3.12 (1) (3)]{corti}$)$. 
\end{itemize}
Then $mD$ is $\pi$-generated for $m\gg 0$, that is, 
there exists a positive integer $m_0$ such that for 
every $m\geq m_0$ the natural homomorphism $\pi^*\pi_*\mathcal O_X(mD)
\to \mathcal O_X(mD)$ is surjective. 
\end{thm}

Before the proof of 
Theorem \ref{thm2}, let us 
recall the definition of {\em{abundant}} divisors, which are called 
{\em{good}} divisors in \cite{kawamata}. See \cite[\S 6-1]{kmm}. 

\begin{defn}[Abundant divisor] 
Let $X$ be a complete normal 
variety and $D$ a $\mathbb Q$-Cartier nef divisor on $X$. 
We define the {\em{numerical Iitaka dimension}} to be 
$$
\nu(X,D)=\max \{e; D^e\not\equiv 0\}. 
$$ 
This means that $D^{e'}\cdot S=0$ for any 
$e'$-dimensional subvarieties $S$ of $X$ with 
$e'>e$ and there exists an $e$-dimensional 
subvariety $T$ of $X$ such that 
$D^e \cdot T>0$. 
Then it is easy to see that $\kappa (X,D)\leq \nu(X, D)$, where 
$\kappa (X,D)$ denotes {\em{Iitaka's $D$-dimension}}. 
A nef $\mathbb Q$-divisor $D$ is said to be {\em{abundant}} if the equality 
$\kappa (X, D)=\nu(X, D)$ holds. 
Let $\pi:X\to S$ be a proper surjective morphism of 
normal varieties and $D$ a $\mathbb Q$-Cartier divisor on $X$. 
Then $D$ is said to be $\pi$-abundant if $D|_{X_\eta}$ is abundant, where 
$X_\eta$ is the generic fiber of $\pi$.   
\end{defn}

\begin{proof}[Proof of {\em{Theorem \ref{thm2}}}] 
If $H-(K_X+B)$ is 
$\pi$-big, then the statement follows from 
the original Kawamata--Shokurov base point 
free theorem. 
Thus, from now on, we assume 
that $H-(K_X+B)$ is not $\pi$-big. 
Then there exists a diagram 
$$
\begin{CD}
Y@>{f}>>Z\\
@V{\mu}VV @VV{\varphi}V\\
X@>>{\pi}>S
\end{CD}
$$
which satisfies the following conditions (see 
\cite[Proposition 6-1-3 and Remark 
6-1-4]{kmm} or 
\cite[Lemma 6]{nakayama}): 
\begin{itemize}
\item[(i)] $\mu, f$ and $\varphi$ are projective morphisms, 
\item[(ii)] $Y$ and $Z$ are non-singular varieties, 
\item[(iii)] $\mu$ is a birational morphism and $f$ is a surjective morphism having 
connected fibers, 
\item[(iv)] there exists a $\varphi$-nef and 
$\varphi$-big $\mathbb Q$-divisor $M_0$ on $Z$ such that 
$$
\mu^*(H-(K_X+B))\sim_{\mathbb Q}f^*M_0, 
$$ 
and 
\item[(v)] there is a $\varphi$-nef $\mathbb Q$-divisor $D$ on $Z$ such that 
$$
\mu^*H\sim _{\mathbb Q}f^*D. 
$$
\end{itemize}
Note that $f:Y\to Z$ is the Iitaka fibration with 
respect to $H-(K_X+B)$ over $S$. 
We put $K_Y+B_Y=\mu^*(K_X+B)$ and $H_Y=\mu^*H$. 
We note that $(Y,B_Y)$ is not necessarily klt but 
{\em{sub}} klt. Thus, 
we have 
$H_Y-(K_Y+B_Y)
\sim_{\mathbb Q}f^*M_0$ (resp.~$H_Y\sim _{\mathbb Q}f^*D$), 
where $M_0$ (resp.~$D$) is a $\varphi$-nef and 
$\varphi$-big (resp.~$\varphi$-nef) $\mathbb Q$-divisor as we 
saw in (iv) and (v). 
Furthermore, we can assume that $D$ and $H$ are Cartier divisors 
and $H_Y\sim f^*D$ 
by replacing $D$ and $H$ by sufficiently divisible 
multiples. 
If we need, we modify $Y$ and $Z$ birationally and 
can assume the following conditions: 
\begin{itemize}
\item[(1)] $K_Y+B_Y\sim _{\mathbb Q}f^*(K_Z+B_Z+M)$, where 
$B_Z$ is the {\em{discriminant}} 
$\mathbb Q$-divisor of $(Y,B_Y)$ on 
$Z$ and $M$ is the {\em{moduli}} 
$\mathbb Q$-divisor on $Z$, 
\item[(2)] $(Z,B_Z)$ is a sub klt pair, 
\item[(3)] $M$ is a 
$\varphi$-nef $\mathbb Q$-divisor on $Z$, 
\item[(4)] $\varphi_*\mathcal O_Z(\ulcorner \mathbf A(Z, B_Z)\urcorner +j
\overline{D})
\subseteq \varphi_*\mathcal O_Z(jD)$ for every 
positive integer $j$, and 
\item[(5)] $D-(K_Z+B_Z)$ is $\varphi$-nef and $\varphi$-big. 
\end{itemize}
Indeed, let $P\subset Z$ be a prime divisor. Let $a_P$ be the 
largest real number $t$ 
such that $(Y,B_Y+tf^*P)$ is sub lc over the generic point of $P$. 
It is obvious that $a_P=1$ for 
all but finitely many prime divisors $P$ of $Z$. 
We note that $a_P$ is a positive rational number for any $P$. 
The {\em{discriminant $\mathbb Q$-divisor}} on $Z$ defined by the following 
formula
$$
B_Z=\sum_P (1-a_P)P. 
$$ 
We note that $\llcorner B_Z\lrcorner \leq 0$. By the properties (iv) and (v), 
we can write 
$$
K_Y+B_Y\sim_{\mathbb Q}f^*(M_1)
$$ 
for a $\mathbb Q$-Cartier 
divisor $M_1$ on $Z$. 
We define $M=M_1-(K_Z+B_Z)$ and call it the {\em{moduli $\mathbb Q$-divisor}} 
on $Z$, where $B_Z$ is the discriminant $\mathbb Q$-divisor defined above. 
Note that $M$ is called the {\em{log-semistable part}} in 
\cite[Section 4]{fujino-mori}.  
So, the condition (1) obviously holds by the definitions of the discriminant 
$\mathbb Q$-divisor $B_Z$ and the moduli $\mathbb Q$-divisor $M$. 
If we take birational modifications of $Y$ and $Z$ suitably, 
we have that $M$ is $\varphi$-nef and $(Z,B_Z)$ is sub klt. 
Thus we obtain (2) and (3). 
For the details, see \cite[Theorems 0.2 and 2.7] {ambro1} or 
Theorem \ref{322} below. 
We note the following lemma 
(cf.~\cite[Lemma 6.2]{ambro1}), which we need to 
apply \cite[Theorems 0.2 and 2.7]{ambro1} or Theorem \ref{322} 
to $f:Y\to Z$ 
(see the condition (2) in \ref{611}). 

\begin{lem}\label{49}
We have $\rank f_*\mathcal O_Y(\ulcorner \mathbf A(Y,B_Y)\urcorner) 
=1$. 
\end{lem}
\begin{proof}[Proof of {\em{Lemma \ref{49}}}]
Since $\mathcal O_Z\simeq 
f_*\mathcal O_Y\subseteq f_*\mathcal O_Y
(\ulcorner \mathbf A(Y, B_Y)\urcorner)$, 
we know $\rank f_*\mathcal O_Y(\ulcorner 
\mathbf A(Y, B_Y)\urcorner )\geq 1$. 
Without loss of generality, we can shrink 
$S$ and assume that $S$ is affine. 
Let $A$ be a $\varphi$-very ample 
divisor such that $f_*\mathcal O_Y(\ulcorner 
\mathbf A(Y, B_Y)\urcorner )\otimes \mathcal O_Z(A)$ is 
$\varphi$-generated. 
Since $M_0$ is a 
$\varphi$-big $\mathbb Q$-divisor on $Z$, we have 
$\mathcal 
O_Z(A)\subset \mathcal O_Z(mM_0)$ for 
a sufficiently divisible positive integer $m$. 
We note that 
$$\pi_*\mu_*\mathcal O_Y(\ulcorner \mathbf A(Y, B_Y)\urcorner 
+\overline{f^*(mM_0)})\simeq \pi_*\mu_*
\mathcal O_Y(f^*(mM_0)), $$ 
where $\overline{f^*(mM_0)}$ is the Cartier 
closure of $f^*(mM_0)$ (see \cite[Example 2.3.12 (1)]{corti}). 
It is because $\mu^*(H-(K_X+B))=H_Y-(K_Y+B_Y)
\sim_{\mathbb Q}f^*M_0$. 
Therefore, 
\begin{eqnarray*}
&& \varphi_*(f_*\mathcal O_Y(\ulcorner \mathbf A(Y,B_Y)\urcorner)
\otimes \mathcal O_Z(A)) 
\\ &\subseteq &
\varphi_*(f_*\mathcal O_Y(\ulcorner \mathbf A(Y, B_Y)\urcorner)
\otimes \mathcal O_Z(mM_0))
\\&\simeq&  \varphi_*\mathcal O_Z
(mM_0). 
\end{eqnarray*}
So, we see that 
$\rank f_*\mathcal O_Y(\ulcorner \mathbf A(Y, B_Y)\urcorner)\leq 
1$. This completes the proof. 
\end{proof}
We know the following lemma 
by Lemma 9.2.2 and Proposition 9.2.3 
in \cite{ambro3} (see also Theorem \ref{322} (a) below). 
\begin{lem}
We have 
$$\mathcal O_Z(\ulcorner \mathbf 
A(Z, B_Z)\urcorner +j\overline D)\subseteq 
f_*\mathcal O_Y(\ulcorner \mathbf A(Y, B_Y)\urcorner 
+j\overline {H_Y})$$ for every integer $j$. 
\end{lem} 
Pushing forward it by $\varphi$, we obtain 
that 
\begin{eqnarray*}
\varphi_*\mathcal O_Z(\ulcorner \mathbf A(Z, B_Z)\urcorner
+j\overline D)& 
\subseteq 
&\varphi_*f_*\mathcal O_Y
(\ulcorner \mathbf A(Y, B_Y)\urcorner +j\overline {H_Y})\\ 
&\simeq &\pi_*\mu_*\mathcal O_Y(\ulcorner \mathbf A(Y, B_Y)\urcorner 
+j\overline {H_Y})\\ 
&\simeq &
\pi_*\mathcal O_X(\ulcorner \mathbf A(X, B)
\urcorner +j\overline {H})\\ 
& \simeq &\pi_*\mathcal O_X(jH)\\
&\simeq &\pi_*\mu_*\mathcal O_Y(jH_Y)\\
&\simeq&\varphi_*f_*\mathcal O_Y(jH_Y)\\
&\simeq &\varphi_*\mathcal O_Z(jD) 
\end{eqnarray*} 
for every integer $j$. 
Thus, we have (4). 
The relation $H_Y-
(K_Y+B_Y)\sim_{\mathbb Q}f^*(D-(K_Z+B_Z+M))$ implies that 
$D-(K_Z+B_Z+M)$ is $\varphi$-nef and $\varphi$-big. 
By (3), $M$ is $\varphi$-nef. 
Therefore, $D-(K_Z+B_Z)=D-(K_Z+B_Z+M)+M$ is $\varphi$-nef and $\varphi$-big. 
It is the condition (5). 
Apply Theorem \ref{thm1} to $D$ on $(Z, B_Z)$. Then we obtain 
that $D$ is $\varphi$-semi-ample. This implies that $H$ is $\pi$-semi-ample. 
This completes the proof. 
\end{proof}

The following corollaries are obvious by Theorem \ref{thm2}. 

\begin{cor}\label{co25}
Let $(X,B)$ be a klt pair and $\pi:X \to S$ a proper 
surjective morphism of normal varieties. 
Assume that $K_X+B$ is $\pi$-nef and $\pi$-abundant. 
Then $K_X+B$ is $\pi$-semi-ample. 
\end{cor}

\begin{cor}\label{co26}
Let $X$ be a complete normal variety such that $K_X\sim _{\mathbb Q}0$. 
Assume that $X$ has only klt singularities. Let $H$ be a nef and 
abundant $\mathbb Q$-Cartier divisor on $X$. 
Then $H$ is semi-ample. 
\end{cor}

We close this section with 
a useful remark. 

\begin{rem}[{cf.~\cite[Remark 3.5]{fujino5}}]
Let $\pi:X\to S$ be a proper surjective morphism 
of normal varieties and 
$D$ a $\pi$-nef and $\pi$-abundant 
Cartier divisor on $X$. 
Then we can easily check that 
$$
\bigoplus _{m\geq 0}\pi_*\mathcal O_X(mD)
$$
is finitely generated if and only if 
$D$ is $\pi$-semi-ample. See, for example, \cite[Lemma 3.10]{fujino5}. 

Let $B$ be an effective $\mathbb Q$-divisor 
on $X$ such that 
$(X, B)$ is klt. 
By \cite{bchm}, we know that 
$$
\bigoplus _{m\geq 0}\pi_*\mathcal O_X(\llcorner 
m(K_X+B)\lrcorner)
$$ 
is finitely generated. 

Assume that $K_X+B$ is $\pi$-nef. 
By the above observation, we obtain that 
$K_X+B$ is $\pi$-semi-ample 
if and only if $K_X+B$ is $\pi$-nef and 
$\pi$-abundant. 
Therefore, we do not need Theorem \ref{thm2} 
to obtain Corollaries \ref{co25} and \ref{co26}. 
\end{rem}

\section{Appendix:~Quick review of Ambro's formula}\label{formula}
In this appendix, we quickly review Ambro's 
formula. For the details, see the original paper \cite{ambro1} 
or Koll\'ar's survey article \cite{kosurvey}. 

\begin{say}\label{611}
Let $f:X\to Y$ be a proper surjective morphism of normal 
varieties and $p:Y\to S$ a proper 
morphism onto a variety $S$. 
Assume the following conditions: 
\begin{itemize}
\item[(1)] $K_X+B$ is 
$\mathbb Q$-Cartier and 
$(X, B)$ is sub klt over the generic point of $Y$, 
\item[(2)] $\rank f_*\mathcal O_X(\ulcorner 
\mathbf A(X, B)\urcorner )=1$, and 
\item[(3)] $K_X+B\sim_{\mathbb Q, f}0$. 
\end{itemize}
By (3), we can write 
$K_X+B\sim _{\mathbb Q}f^*D$ for some $\mathbb Q$-Cartier 
divisor $D$ on $Y$. 
Let $B_Y$ be the {\em{discriminant $\mathbb Q$-divisor}} 
on $Y$. For the definition, see the proof of Theorem \ref{thm2}. 
We put $M_Y=D-(K_Y+B_Y)$ and 
call $M_Y$ the {\em{moduli $\mathbb Q$-divisor}} 
on $Y$. 
Then we have 
$K_X+B\sim_{\mathbb Q}f^*(K_Y+B_Y+M_Y)$. 
Let $\sigma:Y'\to Y$ be a proper birational morphism from a normal 
variety $Y'$. Then we obtain the following commutative diagram: 
$$
\begin{CD}
X@<\mu<< X'\\ 
@V{f}VV @VV{f'}V\\ 
Y@<<{\sigma}<Y' 
\end{CD}
$$
such that 
\begin{itemize}
\item[(i)] $\mu$ is a birational morphism 
from a normal variety $X'$, 
\item[(ii)] we put $K_{X'}+B'=\mu^*(K_X+B)$. Then 
we can write $K_{X'}+B'\sim_{\mathbb Q}{f'}^*(K_{Y'}+B_{Y'}
+M_{Y'})$, where 
$B_{Y'}$ is the discriminant $\mathbb Q$-divisor on $Y'$ associated 
to $f':X'\to Y'$. 
\end{itemize}
Ambro's theorem \cite[Theorems 0.2 and 2.7]{ambro1} says 

\begin{thm}\label{322} 
If we choose $Y'$ appropriately, 
then we have the following properties 
for every proper birational morphism $\nu:Y''\to Y'$ from 
a normal variety $Y''$. 
\begin{itemize}
\item[(a)] $K_{Y'}+B_{Y'}$ is $\mathbb Q$-Cartier and 
$\nu^*(K_{Y'}+B_{Y'})=K_{Y''}+B_{Y''}$. 
In particular, $\mathbf A(Y', B_{Y'})_{Y''}=-B_{Y''}$. 
\item[(b)] The moduli 
$\mathbb Q$-divisor $M_{Y'}$ is nef over 
$S$ and 
$\nu^*(M_{Y'})=M_{Y''}$.   
\end{itemize}
\end{thm}
We note that the nefness of the moduli $\mathbb Q$-divisor 
follows from Fujita--Kawamata's semi-positivity theorem. 
It is a consequence of the theory of variation of Hodge 
structures. For details, see, for example, \cite[Section 5]{mori}, 
\cite[Section 5]{fuji-rem}, or \cite{kosurvey}. 
\end{say}

\section{Kawamata's theorem for varieties in class $\mathcal C$}\label{sec4} 
In this section, we treat Nakayama's 
theorem:~\cite[Theorem 5.5]{nakayama1}, which is 
Kawamata's theorem for varieties in class $\mathcal C$. First, 
let us recall the definition of the varieties in class $\mathcal C$. 

\begin{defn}[{Class $\mathcal C$}]\label{defn1}
A compact complex variety in class $\mathcal C$ is 
a variety which is dominated 
by a compact K\"ahler manifold. 
It is known that $X$ is in class $\mathcal C$ if and 
only if $X$ is bimeromorphically equivalent to a compact 
K\"ahler manifold. 
\end{defn}

Next, we recall the definitions of the K\"ahler cone and 
the nef line bundles on a compact K\"ahler manifold. 

\begin{defn}[K\"ahler cone]\label{defn2} 
Let $Y$ be a $d$-dimensional compact K\"ahler 
manifold. 
We define the K\"ahler cone $\KC(Y)$ of $Y$ to be the set 
$$
\{ [\omega] \in H^{1, 1} (Y, \mathbb R); 
\omega \ {\text{is a K\"ahler form on}} \ Y. \}, 
$$
where $H^{1,1} (Y, \mathbb R):= 
H^2(Y, \mathbb R)\cap H^{1, 1} (Y, \mathbb C)$. 
Then $\KC(Y)$ is an open convex cone in $H^{1, 1} (Y, \mathbb R)$. 
$\overline {\KC}(Y)$ is the closure 
of $\KC(Y)$ in $H^{1, 1} (Y, \mathbb R)$. 
\end{defn}

Finally, we recall the definitions of the quasi-nef line bundles, 
the homological Kodaira dimension, and the quasi-nef and abundant line 
bundles, which were introduced in \cite{nakayama1}. 

\begin{defn}[{cf.~\cite[Definition 2.4]{nakayama1}}]\label{defn3} 
Let $L$ be a line bundle on a compact K\"ahler 
manifold $Y$. $L$ is said to 
be {\em{nef}} if the real 
first Chern class $c_1(L)$ is contained in $\overline{\KC}(Y)$. 
\end{defn}

\begin{rem}
For a new numerical characterization of 
the K\"ahler cone of 
a compact K\"ahler manifold, 
see \cite[Main Theorem 0.1]{dp}. 
A nef line bundle on a compact K\"ahler manifold 
can be characterized numerically by \cite[Corollaries 0.3 and 
0.4]{dp}. 
\end{rem}

\begin{defn}[{cf.~\cite[Definition 2.6]{nakayama1}}]\label{defn4}
Let $X$ be a compact complex variety 
in class $\mathcal C$. 
A line bundle $L$ on $X$ is called 
{\em{quasi-nef}} if there exists a bimeromorphic morphism 
$\mu:Y\to X$ from a compact K\"ahler manifold $Y$ such that 
$\mu^*L$ is nef.  
\end{defn}

\begin{defn}[{cf.~\cite[Definition 2.9]{nakayama1}}]\label{defn5} 
Let $L$ be a quasi-nef line bundle on a complex variety 
$X$ in class $\mathcal C$. Take 
a bimeromorphic morphism $\mu:Y\to X$ from a compact K\"ahler 
manifold $Y$ such that 
$\mu^*L$ is nef. 
Then we define $$\kappa _{\h} (L):= 
\max \{ l \geq 0; 0 \ne c_1(\mu^*L)^l\in 
H^{l, l} (Y, \mathbb R) \}$$ and 
call it the 
{\em{homological Kodaira dimension}} of $L$. 
It is well-defined, because it is independent of the 
choice of $Y$.  
\end{defn}

\begin{defn}[{cf.~\cite[Definition 2.11]{nakayama1}}]\label{defn6} 
Let $L$ be a line bundle 
on a compact complex 
variety $X$ in class $\mathcal C$. $L$ is said to 
be {\em{big}} if $\kappa 
(X, L)=\dim X$. If $L$ is quasi-nef and 
$\kappa (X, L)=\kappa _{\h} (L)$, then $L$ is called {\em{abundant}}. 
\end{defn}

Now, we sate the main theorem of this section. 
It is nothing but \cite[Theorem 5.5]{nakayama1}. 
The reader can find some applications of Theorem \ref{main-da} 
in \cite{cop}. 

\begin{thm}[{cf.~\cite[Theorem 5.5]{nakayama1}}]\label{main-da} 
Let $X$ be a compact normal complex variety in class $\mathcal C$, $B$ an 
effective $\mathbb Q$-divisor 
on $X$, and $H$ a $\mathbb Q$-Cartier divisor 
on $X$. 
Then $H$ is semi-ample under the following conditions: 
\begin{itemize}
\item[(1)] $(X, B)$ is klt, 
\item[(2)] $H$ is quasi-nef, 
\item[(3)] $H-(K_X+B)$ is quasi-nef and abundant, and 
\item[(4)] $\kappa _{\h} (aH-(K_X+B))=\kappa _{\h} 
(H-(K_X+B))$ and $\kappa (X, aH-(K_X+B))\geq 0$ for 
some $a \in \mathbb Q$ with $a>1$. 
\end{itemize} 
\end{thm}

\begin{proof}[Sketch of the proof] 
First, we recall Nakamaya's 
result. 
\begin{lem}[{\cite[Proposition 
2.14 and Corollary 2.16]{nakayama1}}]\label{499} 
There exists the following 
diagram 
$$
\begin{CD}
X@<{\mu}<< Y @>{f}>> Z, 
\end{CD}
$$
where 
\begin{itemize}
\item[(a)] $Y$ is a compact K\"ahler manifold 
and $\mu$ is a bimeromorphic 
morphism, 
\item[(b)] $Z$ is a smooth projective variety, 
\item[(c)] $f$ is surjective and has connected fibers, 
\item[(d)] there exists a nef and 
big $\mathbb Q$-divisor $M_0$ on $Z$ such 
that 
$$
\mu^*(H-(K_X+B))\sim _{\mathbb Q}f^*M_0, 
$$
and 
\item[(e)] there is a nef $\mathbb Q$-divisor 
$D$ on $Z$ such that 
$$
\mu^*H\sim _{\mathbb Q}f^*D. 
$$
\end{itemize}
We note that $Z$ is a smooth {\em{projective}} variety. 
\end{lem}

Let $f:Y\to Z$ be the proper 
surjective morphism from a compact 
K\"ahler manifold $Y$ to a normal projective 
variety $Z$ obtained in Lemma \ref{499}. 
Let $B_Y$ be a $\mathbb Q$-divisor on $Y$ such that 
$K_Y+B_Y=\mu^*(K_X+B)$. 
Then we have the following properties: 
\begin{itemize}
\item[(1)] $K_Y+B_Y$ is 
$\mathbb Q$-Cartier and $(Y_z, B_z)$ is sub klt 
for general $z\in Z$, where 
$Y_z=f^{-1}(z)$ and 
$B_z=B_Y|_{Y_z}$, 
\item[(2)] $\rank f_*\mathcal O_Y(\ulcorner 
\mathbf A (Y, B_Y)\urcorner)=1$, and 
\item[(3)] $K_Y+B_Y\sim _{\mathbb Q, f}0$. 
\end{itemize}
We note that (1) is obvious by the definition of $B_Y$, 
(2) follows from the proof of 
Lemma \ref{49}, and (3) is also 
obvious by Lemma \ref{499}. 
Under these conditions (1), (2), and (3), Ambro's theorem 
(see \cite[Theorems 0.2 and 2.7]{ambro1} or Theorem \ref{322}) holds 
if we use \cite[3.7.~Theorem (4)]{nakayama2} in the 
proof of Ambro's theorem. 
Note that it is not difficult to modify the arguments in \cite{ambro1} 
for our setting. 
More explicitly, let $\sigma: Z'\to Z$ be a proper 
birational morphism 
from a normal projective variety $Z'$. 
If we choose $Z'$ appropriately, 
then we have the following properties for every 
proper birational morphism $\nu:Z''\to Z'$ from 
a normal projective variety $Z''$. 
\begin{itemize}
\item[(a)] $K_{Z'}+B_{Z'}$ is $\mathbb Q$-Cartier 
and $\nu^*(K_{Z'}+B_{Z'})=K_{Z''}+B_{Z''}$, 
where $B_{Z'}$ and $B_{Z''}$ are the discriminant  
$\mathbb Q$-divisors. 
In particular, $\mathbf A(Z', B_{Z'})_{Z''}=-B_{Z''}$. 
\item[(b)] The moduli $\mathbb Q$-divisor $M_{Z'}$ is 
nef and $\nu^*(M_{Z'})=M_{Z''}$.  
\end{itemize} 
For the details and the notation, 
see Section \ref{formula}. 

By applying Ambro's theorem to $f:Y\to Z$, 
the proof of Theorem \ref{thm2} 
works without any modifications. 
We note that $Z$ is a {\em{projective}} variety. 
Thus, we obtain the semi-ampleness of $H$. 
\end{proof} 

\section{Comments for the coming generation}\label{sec5}

The results in \cite{kawamata} 
had already been used in various papers. 
We think that almost all the papers 
only used the main results 
of \cite{kawamata}, that is, 
Theorems 1.1 and 6.1 in \cite{kawamata}. 
Therefore, by this paper, almost all the 
troubles caused by \cite[Theorem 4.3]{kawamata} 
were removed. 
However, 
some authors 
used arguments in \cite{kawamata}. 
We give some comments for the coming generation. 

\begin{say}\label{511} 
As we pointed out in \cite[Remark 3.10.3]{fujino}, 
the proof of \cite[Theorem 4.3]{kawamata} 
is not completed (see also \cite[Theorem 6-1-6]{kmm}). 
We recall the trouble in \cite{kawamata} here for 
the reader's convenience. 

We use the same notation as in the proof of 
Theorem 4.3 in \cite{kawamata}. 
By \cite[Theorem 3.2]{kawamata}, 
${}^{\prime}\!E^{p, q}_1\to {}^{\prime\prime}\!E^{p, q}_1$ are zero 
for 
all $p$ and $q$. 
It does not directly say that 
$$
H^i(X, \mathcal O_X(-\ulcorner L\urcorner))\to 
H^i(D, \mathcal O_D(-\ulcorner L\urcorner))
$$ 
are zero for all $i$. 
So, the proofs of Theorems 4.4, 4.5, 5.1, and 6.1 in \cite{kawamata} 
do not work. 
It is because everything depends on Theorem 4.3 in \cite{kawamata}. 
Thus, we have no rigorous proofs for \cite[Theorems 6-1-8, 
6-1-9]{kmm}. 
In \cite{kawamata}, there seems to be no troubles 
except the proof of Theorem 4.3. 
\end{say}

If someone corrects the proof of 
\cite[Theorem 4.3]{kawamata}, then the following comments 
are unnecessary. 

\begin{say}
In \cite{nakayama}, Nakayama obtained the relative version of 
Kawamata's theorem. The proof 
given there heavily depends on 
Kawamata's original proof. 
So, it does not work by the trouble in \cite[Theorem 4.3]{kawamata}. 
Of course, \cite[Theorem 5]{nakayama} 
is true by our main theorem:~Theorem \ref{thm2}. 
\end{say}

\begin{say}
Section 5 in \cite{nakayama1} contains the same trouble. 
It is because it depends on Kawamata's paper \cite{kawamata}. 
In Section \ref{sec4}, we give a rigorous proof of 
\cite[Theorem 5.5]{nakayama1}. 
\end{say}

\begin{say}
In \cite{fukuda}, Fukuda obtained a slight generalization 
of Kawamata's theorem. 
See \cite[Proposition 3.3]{fukuda}. 
In the final step of the proof of \cite[Proposition 3.3]{fukuda}, 
Fukuda used \cite[Theorem 5.1]{kawamata}. So, Fukuda's original 
proof also has some troubles 
by \cite[Theorem 4.3]{kawamata}. Fortunately, 
we can prove a slight generalization 
of \cite[Proposition 3.3]{fukuda} in \cite[Section 6]{base}. 
\end{say}

%%%%%%%%%%%%%%%%%%%%%%%%%%%%%%%%%
\ifx\undefined\bysame
\newcommand{\bysame|{leavemode\hbox to3em{\hrulefill}\,}
\fi

\end{document}